\input amstex
\documentstyle{amsppt}
\magnification1200
\tolerance=10000
\overfullrule=0pt
\def\n#1{\Bbb #1}

\def\e11{E_{11}}

\def\vk{\varkappa}

\def\be{\beta}

\def\al{\alpha}

\def\la{\lambda}

\def\g{\goth }

\topmatter
\title
Explicit statement of a conjecture on resultantal varieties
\endtitle
\author
S. Ehbauer, A. Grishkov, D. Logachev\footnotemark \footnotetext{E-mail: logachev94{\@}gmail.com \phantom{*********************************************************}}
\endauthor
\NoRunningHeads
\subjclass 11C20, 05A19, 05E99, 11C99, 15A15 \endsubjclass
\abstract The paper [GLZ] "L-functions of Carlitz modules, resultantal varieties and rooted binary trees" is devoted to a description of some resultantal varieties related to L-functions of Carlitz modules. It contains a conjecture that some of these varieties coincide. This conjecture can be formulated in terms of polynomials, namely, in terms of a fact that an explicitly defined polynomial belongs to the radical of the ideal generated by some other polynomials. We give an explicit statement of this conjecture and a numerical result.
\endabstract
\endtopmatter
\document
Paper [GLZ] establishes relations between 3 types of objects: first, L-functions of Carlitz modules, second, resultantal varieties, and third, finite rooted weighted binary trees. The paper is "a garden of conjectures", which look very difficult. 
\medskip
The purpose of the present paper is to attract the attention of the mathematical community to the first (in logical order) of these conjectures ([GLZ], Conjecture 0.2.4 = [GL], Conjecture 9.3). Here we give an elementary (as explicit as possible) statement of a particular case of this conjecture. Also, we give numerical examples. 
\medskip
This particular case is the following. Conjecture 0.2.4 of [GLZ] depends on an integer parameter $n>0$; we state the conjecture only for the case $n=1$. For $n>1$ see Remark 14. 
\medskip
The reason is the following. The statement of Conjecture 0.2.4 in [GLZ] is overwhelmed in details and is hardly understandable for a beginner. It states that a polynomial belongs to the radical of an ideal in a ring of polynomials. So, here we give an explicit definition of the polynomial and the ideal.  See [GLZ] for details, significance of the problem and possible generalizations. We mention that --- because both the polynomial and the ideal --- come from determinants, a good reference is [K].  
\medskip
Let $m\ge2$ be integer and $a_*=(a_0, ...,a_m)$ independent variables, $a_\vk=0$ for $\vk\not\in \{0,\dots,m\}$. The matrix $\g M(m)(a_0,...,a_m)$ is a $(m-1)\times (m-1)$-matrix whose $(\al, \be)$-th entry is equal to $a_{2\be-\al}$. If it is clear what objects $(a_0, ...,a_m)$ are kept in mind, we write $\g M(m)$ instead of $\g M(m)(a_0,...,a_m)$.
\medskip
{\bf Example.} For $m=4, \ 5, \ 6, \ 7$ \ \ $\g M(m)$ are the following (their structure is slightly different for odd and even $m$):
\medskip
\newpage
\centerline{$\g M(4)(a_0,...,a_4)$  \ \ \ \ \ \ \ \ \ \ \ \ \ \ \ \ \ \ \ \  $\g M(5)(b_0,...,b_5)$ }
\medskip
$$\left(\matrix a_1&a_3&0\\ a_0&a_2&a_4 \\ 0&a_1&a_3\endmatrix \right)
\ \ \ \ \ \  \ \ \ \ \ \ \ \ \ \ \ \ \ \ \ \left(\matrix b_1&b_3&b_5&0\\ b_0&b_2&b_4&0\\ 0&b_1&b_3&b_5 \\ 0&b_0&b_2&b_4\endmatrix \right) \eqno{(1)}$$
\medskip
\centerline{$\g M(6)(a_0,...,a_6)$  \ \ \ \ \ \ \ \ \ \ \ \ \ \ \ \ \ \ \ \ \ \ \ \ $\g M(7)(b_0,...,b_7)$ }
\medskip
$$\left(\matrix a_1&a_3&a_5&0&0\\ a_0&a_2&a_4&a_6&0 \\ 0&a_1&a_3&a_5&0 \\ 0&a_0&a_2&a_4&a_6\\ 0&0&a_1&a_3&a_5\endmatrix \right)
\ \ \ \ \ \   \left(\matrix b_1&b_3&b_5&b_7&0&0\\ b_0&b_2&b_4&b_6&0& 0 \\ 0&b_1&b_3&b_5&b_7&0 \\ 0&b_0&b_2&b_4&b_6&0\\ 0&0&b_1&b_3&b_5&b_7\\ 0&0&b_0&b_2&b_4&b_6 \endmatrix \right) \eqno{(1a)}$$
\medskip
We see that $\g M(m)$ is obtained by a permutation of lines of the Sylvester matrix (defining the resultant) of polynomials whose coefficients are $a_1, \ a_3, \ a_5,\dots$ and $a_2, \ a_4, \ a_6,\dots$.
\medskip
Further, $\g M(m)$ is obtained by the symmetry with respect to the antidiagonal from the Schur matrix for the case $\la_1=m$, $\la_2=m-1, \dots$, see (1.5) of [FP]. 
\medskip
{\bf Remark 1b.} We consider in [GL] a more general situation: the situation of the present paper corresponds to $q=2$, while in [GL] we have: $q$ can be any power of a prime. An analog of $\g M(m)$ for any $q$ (see [GL], (3.2) for an analog of the below $\Cal M(m)$) is the symmetry with respect to the center from the Schur matrix for the case $\la_1=q-1+c$, $\la_2=2(q-1)+c, \ \la_3=3(q-1)+c, \dots$ where $c$ is a constant. An analog of Conjecture 9 for this case exists. 
\medskip
Let $Ch(\g M(m))$ be the $(-1)^{m-1}\cdot$ characteristic polynomial of $\g M(m)$:
$$Ch(\g M(m))=|\g M(m)-U\cdot I_{m-1}|$$
\medskip
{\bf Definition 2.} $D(m,i)\in \n Z[a_0,\dots, a_m]$ are coefficients at $U^i$ of $Ch(\g M(m))$ considered as a polynomial in $U$:
$$Ch(\g M(m))=D(m,0)+D(m,1)\ U+D(m,2)\ U^2+\dots + D(m,m-2)\ U^{m-2}+(-U)^{m-1}\eqno{(3)}$$
They are homogeneous polynomials of degree $m-1-i$. If we define a weight function on $a_0, \dots, a_m$ as follows: $wt(a_i)=i$ then all terms of $D(m,0)$ have the same weight $m(m-1)/2$. 
\medskip
Particularly, we have: $D(m,0)=|\g M(m)|$;
\medskip
$D(m,m-2)=(-1)^m tr (\g M(m))=(-1)^m(a_1+a_2+...+a_{m-1})$;
\medskip
$D(m,m-1)=(-1)^{m-1}$.
\medskip
\newpage
{\bf Example 4.} For $m=4$ we have
\medskip
(a) $D(4,0)=-a_0a_3^2-a_1^2a_4+a_1a_2a_3$;
\medskip
(b) $D(4,1)=a_0a_3-a_1a_2-a_1a_3+a_1a_4-a_2a_3$;
\medskip
(c) $D(4,2)=a_1+a_2+a_3$.
\medskip
Let $t$ be an independent variable. The matrix $\Cal M(m)(a_0,...,a_m;t)$ is a $m\times m$-matrix whose $(\al, \be)$-th entry is equal to 
$$a_{2\be-\al}t-a_{2\be-\al-1}$$ (this is [GLZ], (0.2.0) for $q=2$, $n=1$, $\vk=m$, $P=(a_0, ...,a_m)$ ). For $m=4, \ 5$ we have respectively 
$$\Cal M(4)=\left(\matrix a_1t-a_0&a_3t-a_2&-a_4&0\\ a_0t&a_2t-a_1&a_4t-a_3&0 \\ 0& a_1t-a_0&a_3t-a_2&-a_4\\ 0 &  a_0t&a_2t-a_1&a_4t-a_3\endmatrix \right)$$
$$\Cal M(5)=\left(\matrix a_1t-a_0&a_3t-a_2&a_5t-a_4&0&0\\ a_0t&a_2t-a_1&a_4t-a_3&-a_5&0 \\ 0& a_1t-a_0&a_3t-a_2&a_5t-a_4&0\\ 0 &  a_0t&a_2t-a_1&a_4t-a_3&-a_5\\ 0&0& a_1t-a_0&a_3t-a_2&a_5t-a_4\endmatrix \right)$$
We consider its characteristic polynomial in variable $T$:
$$Ch(\Cal M(m),T):=\det(I_{m} - \Cal M(m)T)\in\n Z[a_0,\dots, a_m][t][T]\eqno{(5)}$$
(this is [GLZ], (0.2.1) for $q=2$, $n=1$, $\vk=m$, $P=(a_0, ...,a_m)$ ). 
\medskip
The polynomials $$H_{ij}=H_{ij}(m)=H_{ij}(a_0,\dots,a_m)\in\n Z[a_0,\dots, a_m]$$ are coefficients of $Ch(\Cal M(m),T)$ at $t^jT^{m-i}$:

$$Ch(\Cal M(m),T)=\sum_{i=0}^m \sum_{ j=0}^{m-i}H_{ij}(m)t^ jT^{m-i}\eqno{(6)}$$

Particularly, we have: $H_{m0}(m)=1$
\medskip
$H_{m-1,0}(m)=-(a_1+...+a_m)$
\medskip
$H_{m-1,1}(m)=a_0+...+a_{m-1}$
\medskip
$H_{0j}(m)=\pm a_j \ D(m,0)$ (see 10.1 below).
\medskip
$H_{ij}(m)$ are homogeneous polynomials in $a_0,...,a_m$ of degree $m-i$. 
\medskip
{\bf 7.} There is a symmetry between $H_{ij}$ and $H_{i,m-i-j}$: the substitution $a_k \mapsto a_{m-k}$ moves $H_{ij}$ to $H_{i,m-i-j}$ (up to a sign). 
\medskip
Hence, for a fixed $m$ the set of all $H_{ij}$ is a triange. For $m=4$ it is the following: 
\medskip
\ \ \ Degree
$$\matrix 0&\ \ &&&&&&H_{40}=1&&&&\\ \\ 1&\ \ &&&&&H_{30}&&H_{31}&&& \\ \\ 2&\ \ &&&&H_{20}&&H_{21}&&H_{22}&& \\ \\ 3&\ \ &&&H_{10}&&H_{11}&&H_{12}&&H_{13}& \\ \\ 4&\ \ &&H_{00}&&H_{01}&&H_{02}&&H_{03}&&H_{04}\endmatrix \eqno{(8)}$$
\medskip
The symmetry (7) is the symmetry with respect to the vertical axis of the triangle (8). 
\medskip
Now we can formulate the Conjecture 0.2.4 of [GLZ] (case $n=1$):
\medskip
{\bf Conjecture 9.} $\forall \ m,\ i\in [0,\dots, m-1]$, $ j\in [0,\dots, m-i]$) $$\exists \ \varkappa \hbox{ such that }(H_{ij}(m))^\varkappa\in \ <D(m,0),\dots,D(m,i)>\eqno{(10)}$$ --- the ideal generated by
$D(m,0),\dots,D(m,i)$.
\medskip
The conjecture is proved for $H_{ij}(m)$ on the sides of the triangle (8), any $m$, see 10.1 - 10.3, and for a point near the upper vertex, see 10.4. Let us give these results (the signs in the below formulas depend on parities of entries, we do not give here their values): 
\medskip
{\bf 10.1.} If $i=0$ then $H_{0j}(m)=\pm a_j \ D(m,0)$ (very simple; see [GL], Theorem III);
\medskip
{\bf 10.2.} If $j=0$ then $H_{i0}(m)=\pm a_0 \ D(m,i)\pm D(m, i-1)$ (also very simple; see [GL], Proposition 9.14);
\medskip
{\bf 10.3.} If $j=m-i$ then $H_{i, m-i}(m)=\pm a_m \ D(m,i)\pm D(m, i-1)$ (follows from (10.2) by the symmetry (7));
\medskip
{\bf 10.4.} $H_{m-2,1}(m) = D(m, m-3)- D(m, m -2)^2$ (requires a few calculations, see [ELS], Theorem 2.1). 
\medskip
We see that for all these cases we have $\vk=1$. Conjecturally, for $n=1$ these are the only cases where $\vk=1$. There is also one more case for $n=2$, see Remark 14. 
\medskip
For the reader's convenience, we reproduce here the table of the minimal $\vk$ as a function of $m$, $i$ (within the limits of the table $\vk$ does not depend on $j$ except $j=0, \ m-i$ where $\vk=1$), see [GLZ], A1 for the original. Recall that for $i=m-2$ we have $\vk=1$. This table was obtained by using the Magma computer algebra system, because it has an option to check whether a polynomial belongs to an ideal, or not. 
\medskip
\newpage
\centerline{\bf Table of $\vk$ as a function on $m$, $i$}
\settabs 8 \columns
\medskip
\medskip
\medskip
\+ $m$ &$i$& & $\vk$\cr
\medskip
\+4&1&    &2\cr
\medskip
\+5&1&&2\cr
\medskip
\+5&2&   &3\cr
\medskip
\+6&1&&2\cr
\medskip
\+6&2&&4\cr
\medskip
\+6&3&  &6\cr
\medskip
\+7&1&&2\cr
\medskip
\+7&2&&4\cr
\medskip
\+7&3&  &? \ ($\ge6$)\cr
\medskip
\+7&4&  &? \cr
\medskip
We could not finish calculations for the last two cases, because of a lack of computer power. 
\medskip
This table gives us 
\medskip
{\bf Conjecture 11.} For $H_{ij}$ of the interior of the triangle (8), any $m$, we have: $\vk$ does not depend on $j$, but only on $i$; for $i=1$ we have $\vk=2$ (clearly here we mean the minimal possible value of $\vk$ from (10)). 
\medskip
{\bf Remark.} We see that for a fixed $m$ we have: $\vk$ grows while $i$ grows from 0 to $m-3$. Nevertheless, for $i=m-2$ we have $\vk=1$, see 10.4. 
\medskip
We get that the first non-trivial case (i.e. $\vk>1$), as well as the only non-trivial case for $m=4$, is the case of $H_{11}(4)$ or, equivalently, $H_{12}(4)$. We calculated:
\medskip
{\bf 12. Result.} $H_{12}(4)^2= (C_1 +a_3 C_2)D(4,1) - C_2 D(4,0)$ where 
\medskip
$C_1=2a_0a_1a_2a_4+a_0a_1a_3a_4-a_0a_1a_4^2+a_0a_2^2a_3+2a_0a_2a_3^2-3a_0a_2a_3a_4+a_0a_3^3$
\medskip
$-2a_0a_3^2a_4 +a_0a_3a_4^2-a_1^2a_2a_4-a_1^2a_3a_4-a_1a_2^3-a_1a_2^2a_3-a_1a_2a_3^2+a_1a_2a_3a_4$
\medskip
$-a_1a_3^3+a_1a_3^2a_4+a_2^3a_3 -a_2a_3^3$
\medskip
$C_2=a_0^2a_4-2a_0a_1a_4-3a_0a_2a_4+a_0a_4^2+a_1^2a_4+a_1a_2^2+a_1a_2a_4+a_2^3$
\medskip
Recall that all involved polynomials are homogeneous, of degrees:
\medskip
deg $H_{12}(4)=3$; deg $D(4,1)=2$; deg $D(4,0)=3$; deg $C_1=4$; deg $C_2=3$. 
\medskip
{\bf 13. Remarks.} 1. $C_2$ does not contain $a_3$. 
\medskip
2. We can specialize the calculations: $a_3 \mapsto 0$. After such specialization, $C_1$ contains only 4 terms. 
\medskip
In order to prove the conjecture, it is desirable to know beforehand the coefficients $\la_k$ of a linear combination 
$$(H_{ij}(m))^\varkappa=\sum_{k=0}^i \la_k \ D(m,k)$$
They are not unique; we are not sure that the above $C_1, \ C_2$ are the simplest ones. Most likely $\la_k$ come from some determinants; but from which ones? 
\medskip
{\bf 14. Remark.} For the statement of the conjecture for $n>1$ see [GLZ], (0.2.4). We indicate here that we must consider an analog of the matrix $\Cal M(m)(a_0,...,a_m;t)$; its entries are polynomials in $t$ of degree $n$ (see [GLZ], (0.1.4a) and (0.2.0)). We have 
\medskip
{\bf Theorem 14.1.} ([GLZ], Proposition 7.2.1). If Conjecture 9 holds for $n=1$ for all $i,\ j,\ m$, then its analog for any $n$ also holds. 
\medskip
Since at the moment we do not know a proof of Conjecture 9, we are interested to prove particular cases of Conjecture for $n>1$. For $n=2$ there is a case where $\vk=1$, namely, $H_{11}$ is a linear combination of $D(m,0)$ and $D(m,1)$, see [ELS], Theorem 3.1 (the proof is 8 pages long). 
\medskip
{\bf References}
\medskip

[ELS] S. Ehbauer, D. Logachev, M. Sarraff de Nascimento. Some cases of a conjecture on L-functions of twisted Carlitz modules. Comm. Algebra 46 (2018), no. 5, 2130 -- 2145. https://arxiv.org/pdf/1707.04339.pdf
\medskip
[FP] Fulton, W., Pragacz, P. (1998). Schubert Varieties and Degeneracy Loci, Lecture Notes in Mathematics, 1689. Berlin: Springer-Verlag.
\medskip
[GL] A. Grishkov, D. Logachev. Resultantal varieties related to zeroes of L-functions of Carlitz modules. Finite Fields and Their Applications. 2016, vol. 38, p. 116 -- 176. https://arxiv.org/pdf/1205.2900.pdf
\medskip
[GLZ] A. Grishkov, D. Logachev, A. Zobnin. L-functions of Carlitz modules, resultantal varieties and rooted binary trees - I.  J. of Number Theory,  2022, vol. 238, p. 269 -- 312. https://arxiv.org/pdf/1607.06147.pdf (a complete version). 
\medskip
[K] Krattenthaler, Christian. Advanced determinant calculus. 

arxiv.org/pdf/math/9902004v3.pdf
\enddocument